\documentclass{amsart}
\usepackage[all]{xy}
\newtheorem{thm}{Theorem}[section]

\newtheorem{prob}{Problem}
\theoremstyle{definition}

\theoremstyle{remark}

\numberwithin{equation}{section}
\begin{document}

\title{On  torsion images of Coxeter groups and  question of Wiegold}
\author{Rostislav Grigorchuk}
\address{Department of Mathematics
Mailstop 3368 Texas A\&M University College Station, TX 77843-3368,
USA }

\email{grigorch@math.tamu.edu}

%\thanks{The author expresses his thanks to T.~Januszkiewicz, V.~Nekrashevych, Z.~{\v S}uni\'c, and I.~Subbotin  for their valuable remarks}
.

%\date{}%
%\dedicatory{Dedicated to Leonid Kurdachenko on the occasion of his  60th birthday.}%
%\commby{}%
% ----------------------------------------------------------------
\begin{abstract}
We show that every Coxeter group that is not virtually abelian and for which all labels in the corresponding
Coxeter graph are powers of 2 or infinity can be mapped onto uncountably many infinite 2-groups which, in
addition, may be chosen to be just-infinite, branch groups of intermediate growth. Also we answer affirmatively a question raised by Wiegold in Kourovka Notebook.
\end{abstract}
\maketitle
% ----------------------------------------------------------------

\section{Introduction}
This note is a shortened and modified version of the publication \cite{grigorch:wiegold} (there is a free access to the content of the journal ``Algebra and Discrete Mathematics'' via
http://adm.lnpu.edu.ua/index.htm ).

One of the most outstanding problems in Algebra known as the Burnside Problem (on periodic groups) was
formulated by Burnside in 1902 and was later split into three branches: the  General Burnside Problem, the
Bounded Burnside Problem, and the Restricted Burnside Problem.  The General Burnside Problem was asking if
there exists an infinite finitely generated torsion group. It was answered positively by Golod in
1964~\cite{golod:p-groups} based on Golod-Shafarevich Theorem~\cite{golod_s:class_field_tower}.  The Bounded
Burnside Problem was solved by S.~P.~Novikov and S.~I.~Adjan~\cite{novikov_ad:burnside2,adian:b-burnside}.
The Restricted Burnside Problem was solved by
E.~Zelmanov~\cite{zelmanov:burnside_odd,zelmanov:burnside_2-groups} as a corollary of his fundamental results
on Lie and Jordan algebras.  The  problem of Burnside inspired a lot of activity and new directions of
research. For solution of these problems, various constructions, and surveys we
recommend~\cite{adian:b-burnside,aleshin:burnside,sushch:burnside,olshanskii:noetherian,grigorch:burnside,grigorch:degrees,gupta_s:burnside,%
kostrikin:burnside,zelmanov:burnside_odd,zelmanov:burnside_2-groups,ivanov:burnside94,lysenok:burneven,gupta:finiteorder,%
vaughan-lee:restricted,ivanov_ol:97,bartholdi_s:growth,grigorchuk-l:burnside,bar_gs:branch,olshanskii_Os:large,
delzant_Grom:burnside,osin:rankgrad} which contain  further information on this topic.

Among various problems around the Burnside problem is the problem on minimal values of periods of elements. In the case of the  Bounded Burnside Problem the main remaining open question is: what is the minimal $n$ such that the free Burnside group
 \[B(m,n)=\langle a_1,\dots, a_m \|\ X^n=1 \rangle\]
given by $m \geq 2 $ generators and the identity $X^n=1$ is infinite? Is it 5,7,8 or a larger number? (it is
known that the exponents 3,4, and 6 produce finite free Burnside groups). By the  celebrated result of
E.~Zelmanov~\cite{zelmanov:burnside_odd,zelmanov:burnside_2-groups} finitely generated torsion group with
bounded periods of elements cannot be residually finite. Therefore in a finitely generated residually finite
torsion group  periods of elements  are not uniformly bounded and  one can study the growth of the period
function as  was initiated in~\cite{grigorch:degrees}. For instance, the group $\mathcal{G}$, constructed by
the author in~\cite{grigorch:burnside} as a simple example of a residually finite 2-group, has polynomial
growth of periods and is just-infinite (i.e. it is infinite but every proper quotient of it is finite).
Therefore making the order of any element of $\mathcal{G}$ smaller will make the group finite.

Fixing the number of generators $m \geq 2$ one may be interested in the minimal values of orders of
generators, of products of their powers, of products of length 3 etc, that $m$-generated infinite residually
finite torsion group may have. The case of  $p$-groups is of special interest because of many reasons. For
$m=2$ the order 2 for the generators $x,y$ is impossible because the group would be a dihedral group in this
case. As we will see, the orders 2 and 4 (and 8 for the product $xy$) are possible values, while the triple
$2,4,4$ is not possible (because the corresponding group is crystallographic). Starting with $m=3$ the orders
of generators may take the minimal possible value 2, and we come to the question on torsion quotients of
Coxeter groups, which is the main topic of this note. As Coxeter groups are generated by involutions it is
natural to investigate their 2-torsion quotients.

Recall that a Coxeter group can be defined as a group with a presentation
\[\mathcal{C}=\langle x_1,x_2,\dots , x_n \|\ x_i^2, (x_ix_j)^{m_{ij}}, 1
\leq i < j \leq n \rangle, \] where $m_{i,j} \in {\mathbb N} \cup \{\infty\}$ (the case $m_{i,j}=\infty$
means that there is no defining relator involving $x_i$ and $x_j$).

If $m_{ij}=2$ this means that $x_i$ and $x_j$ commute. A Coxeter group can be described by a Coxeter graph $\mathcal{Z}$. The vertices of the graph are labeled by the generators of the group $\mathcal{C}$, the vertices $x_i$ and $x_j$ are connected by an edge if and only if $m_{i,j} \geq 3$, and an edge is labeled by the corresponding value $m_{ij}$ whenever this value is 4 or greater. If a Coxeter graph is not connected, then the group $\mathcal{C}$ is a direct product of Coxeter subgroups corresponding to the connected components. Therefore we may focus on the case of connected Coxeter graphs. If we are  interested in 2-torsion quotients of $\mathcal{C}$, then one has to assume that $m_{ij}$ are powers of 2 or infinity. In order for $\mathcal{C}$ to have infinite torsion quotients it has to be infinite and not virtually abelian. The list of finite and virtually abelian Coxeter groups with connected Coxeter graphs is well known. A
comprehensive treatment of Coxeter groups can be found in M.~Davis' book \cite{davis:coxeter}.

\begin{thm}\label{coxeter} Let $\mathcal{C}$ be a non virtually abelian
Coxeter group defined by a connected Coxeter graph $\mathcal{Z}$ with all edge labels $m_{ij}$ being powers of 2 or infinity. If  $\mathcal{Z}$ is not a tree or is a tree with $\geq 4$ vertices, or is a tree with two edges with one label $\geq 4$ and the other $\geq 8$, then the group $\mathcal{C}$ has uncountably many 2-torsion quotients. Moreover these quotients can be chosen to be residually finite, just-infinite, branch 2-groups of intermediate growth and the main property that distinguishes them is the growth type of the group.
\end{thm}

Observe that all cases of connected Coxeter graphs that are excluded by the statement of
Theorem~\ref{coxeter} are related to finite or virtually abelian crystallographic groups. Indeed, in the case when $\mathcal{Z}$ consist of one edge the corresponding group is a dihedral group, and when $\mathcal{Z}$ has two edges labeled by 4 the corresponding Coxeter group is  the crystallographic group $\langle x,y,z \|\
x^2,y^2,z^2,(yz)^2,(xy)^4,(xz)^4 \rangle$ generated by reflections in sides of an isosceles right triangle.

On the other hand, there are four ``critical'' Coxeter groups $\Xi$, $\Phi$, $\Upsilon$, and $\Pi$:
\[\Xi= \langle a,c,d \|\ a^2,c^2,d^2,(cd)^2,(ad)^4,(ac)^8 \rangle,\]
\[\Phi=\langle x,y,z\|\ x^2,y^2,z^2,(xy)^4,(xz)^4,(yz)^4 \rangle,\]
\[\Upsilon=\langle a,b,c,d \|\
a^2,b^2,c^2,d^2,(ac)^2,(ad)^2,(bd)^2,(ab)^4,(bc)^4,(cd)^4 \rangle,\]
\[\Pi=\langle a,b,c,d \|\
a^2,b^2,c^2,d^2,(bc)^2,(bd)^2,(cd)^2,(ab)^4,(ac)^4,(ad)^4 \rangle,\]
that satisfy the requirements of Theorem~\ref{coxeter} and play a crucial role in the proof.
Their Coxeter graphs are depicted in Figure~1.

\begin{figure}[!ht]
\[
\xymatrix@C=25pt@R=20pt{
\Xi&&&&&\Phi& *[o][F-]{\bullet} \ar@{-}[ddl]_{4} \ar@{-}[ddr]^{4} \ar@{}[r]|<<<{x}& \\
*[o][F-]{\bullet} \ar@{-}[r]^{4} \ar@{}[d]|<<<{c}& *[o][F-]{\bullet} \ar@{-}[r]^{8} \ar@{}[d]|<<<{a} & *[o][F-]{\bullet} \ar@{}[d]|<<<{d} &&&&&\\
&&&&&*[o][F-]{\bullet} \ar@{-}[rr]_{4} \ar@{}[u]|<<<{y}&& *[o][F-]{\bullet} \ar@{}[u]|<<<{z}\\
&&&&&*[o][F-]{\bullet} \ar@{-}[dr]^{4} \ar@{}[d]|<<<{b}&& *[o][F-]{\bullet} \ar@{-}[dl]_{4} \ar@{}[d]|<<<{c}\\
*[o][F-]{\bullet} \ar@{-}[r]^{4} \ar@{}[d]|<<<{a}& *[o][F-]{\bullet} \ar@{-}[r]^{4} \ar@{}[d]|<<<{b} & *[o][F-]{\bullet} \ar@{-}[r]^{4} \ar@{}[d]|<<<{c} & *[o][F-]{\bullet} \ar@{}[d]|<<<{d}&&&*[o][F-]{\bullet} \ar@{-}[d]_{4} \ar@{}[dr]|<<<{a}&  \\
\Upsilon&&&&&\Pi&*[o][F-]{\bullet} \ar@{}[r]|<<<{d}&
 }
\]
\caption{Coxeter graphs corresponding to $\Xi$, $\Phi$, $\Upsilon$, and $\Pi$}
\end{figure}
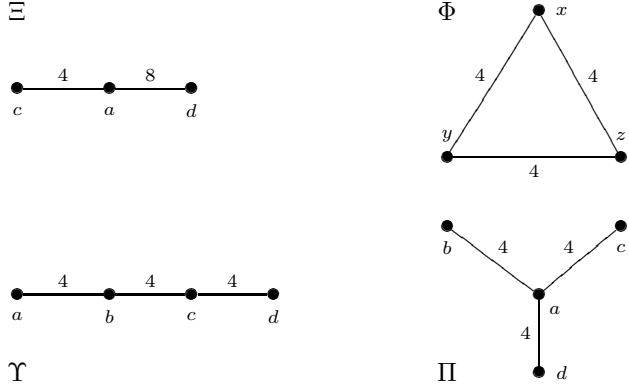

The proof of the theorem is based on the properties of the group $\mathcal{G}$ and of the groups of
intermediate growth from the uncountable family $\{G_{\omega} \mid \omega \in \Omega\}$ constructed in
\cite{grigorch:degrees}, which includes (and generalizes) the example $\mathcal{G}$ (some information about
groups $G_{\omega}$ will be provided below).

The definition of a branch group is a bit involved and we direct the reader to
\cite{grigorch:jibranch,grigorch:branch,bar_gs:branch} for more information on branch groups. A group $G$ is
a branch group if it has a strictly decreasing sequence $\{H_n\}_{n=0}^{\infty}$ of normal subgroups of
finite index with trivial intersection, satisfying the following properties:
\[ [H_{n-1}:H_n]=m_n \in \mathbb{N},\]
for $n=1,2,\dots$, there is a decompositions of $H_n$ into the direct product of $N_n=m_1m_2\dots m_n$ copies
of a group $L_n$ such that the decomposition for $H_{n+1}$ refines the decomposition for $H_n$  (in the sense
that each factor of $H_n$ contains the product of $m_{n+1}$ factors of the decomposition of $H_{n+1}$), and
for each $n$ the group $G$ acts transitively by conjugation on the set of factors of $H_n$. Branch groups
constitute one of three classes into which the class of just-infinite groups naturally splits and they appear
in various situations \cite{grigorch:jibranch, bartholdi_gn:fractal,
nekrash:self-similar,bartholdi_nekr:quandratic}.

The natural language to work with branch groups is via their actions on regular rooted trees as described in
\cite{grigorch:jibranch,grigorchuk-s:standrews,bar_gs:branch}. Then, by definition, a group $G$ acting by
automorphisms on a binary rooted tree $T$ (without change the definition holds also for arbitrary spherically
homogenous rooted tree) is branch if it acts transitively on levels and for any $n \geq 1$ the rigid
stabilizer $rist_G(n)$ of level $n$ has finite index in $G$. (Rigid stabilizer of level $n$  is the subgroup
generated by the rigid stabilizers of the vertices at level $n$, and the rigid stabilizer of a vertex $u$  is
the subgroup of $G$ acting trivially outside the subtree $T_u$ with root $u$). Observe that $rist_G(n)$ is
the direct product of $rist_G(u)$, where $u$ runs over the set of vertices of level $n$, which makes a link
to the algebraic definition given before.

Let $st_G(1)$ be the stabilizer of the first level. Then $\psi : st_G(1) \rightarrow A\times B$ is an embedding, where $A$ and $B$ are the projections of $G$ on the left and right, respectively, rooted subtree of $T$ with roots at the first level.

The groups $G_\omega, \omega \in \Omega_1$ (the sets $\Omega,\Omega_0,\Omega_1$ will be defined later), and in particular the group $\mathcal{G}$, are branch, just-infinite
groups~\cite{grigorch:degrees,grigorch:jibranch} (the term \emph{branch group} is not used in
\cite{grigorch:degrees} as at the time of writing of the paper there was no definition of this class of
groups, but the proof of \cite[Theorem~2.2]{grigorch:degrees}  implies the branch property).

The subgroups $\langle b,ac \rangle, \langle c,ad \rangle,\langle d,ab \rangle$  of index 2 in $\mathcal{G}$,
and the corresponding subgroups $\langle b_{\omega},ac_{\omega} \rangle, \langle c_{\omega},ad_{\omega}
\rangle,\langle d_{\omega},ab_{\omega} \rangle$  of index 2 in $G_\omega, \omega \in \Omega_0$ are also
branch, because they act transitively on binary tree $T$ as one can easy check or apply the criterion
from~\cite[Theorem~2]{grigorch:jibranch}. As any proper quotient of a branch group is virtually
abelian~\cite[Theorem~4]{grigorch:jibranch}, and as all  groups  $G_\omega, \omega \in \Omega_0$ are branch
2-groups, they are just-infinite, as well as are just-infinite the subgroups of index 2 listed above.

A finitely generated group has intermediate growth if the growth function $\gamma(n)$, counting the number of elements of length at most $n$, grows faster than any polynomial but slower than any exponential function $\lambda^n$, for $\lambda > 1$. We use Milnor's equivalence on the set of growth functions of finitely
generated groups: $\gamma_1(n)\sim \gamma_2(n)$ if there is $C \in \mathbb{N}$ such that $\gamma_1(n) \leq
\gamma_2(Cn)$ and $\gamma_2(n)\leq \gamma_1(Cn)$, for $n=0,1,2,\dots$. For a given finitely generated group the class of equivalence of its growth function does not depend on the choice of a finite generating set and is called the growth degree of the group.  It is shown in~\cite{grigorch:degrees} that there are uncountably many growth degrees of finitely generated groups and, moreover, the partially ordered set of growth degrees of finitely generated groups contains both chains and antichains of continuum cardinality. Some additional information about the growth properties of the family $\{G_\omega\}, \omega \in \Omega_1$ will be provided in
the next section.

\section{Preliminary facts}

The group $\mathcal{G}$  was defined in~\cite{grigorch:burnside} as a group generated by four interval exchange transformations $a,b,c,d$ of order 2 acting on the interval $[0,1]$ from which the diadic rational points are removed. From the definition it immediately follows that the generators satisfy the relations
\[a^2=b^2=c^2=d^2=[b,c]=[b,d]=[c,d]=bcd=(ad)^4=(ac)^8=(ab)^{16}=1\]
(this list of relations is not complete).
The branch algorithm for decision of the word problem
described in \cite{grigorch:degrees} is very efficient and has time (or space) complexity $n\log(n)$.
As shown by I.~Lys{\"e}nok~\cite{lysionok:presentation}, $\mathcal{G}$ can be described by the following
presentation
\begin{equation} \label{Lys}
\langle a,b,c,d\
| a^2,b^2,c^2,d^2, bcd, \alpha^n((ad)^4), \alpha^n((adacac)^4), \ n \geq 0 \rangle,
\end{equation}
where $\alpha$ is the substitution $\alpha: a\rightarrow aca, b \rightarrow d, c \rightarrow b, d \rightarrow
c$. It is very interesting and surprising that the relators in Lys\"{o}nok presentations are words of power
at most 8, as we know nothing about the free Burnside group of exponent 8. The group $\mathcal{G}$ is not
finitely presented, and it is shown in~\cite{grigorchuk_schur} that the relators given in ~\eqref{Lys} are
independent (i.e., none of them can be deleted from the set of relators without changing the group). The
relation $bcd=1$ implies that the group $\mathcal{G}$ is 3-generated, but it is usually convenient to work
with the generating set $\mathcal{A}=\{a,b,c,d\}$, because together with the identity element it constitutes
the so called nucleus of the group, an important tool in the study of self-similar groups
\cite{nekrash:self-similar}. Excluding the generator $b$ we see that $\mathcal{G}$ is a homomorphic image of the group $\Xi$.

 For the proof of Theorem~\ref{coxeter} we will use the construction of
an uncountable family of groups $G_\omega$, where  $\omega \in \Omega=\{0,1,2\}^{\textbf{N}}$ described in
\cite{grigorch:degrees} for which the group $\mathcal{G}$ is a particular case corresponding to the sequence
$\varsigma=(012)^{\infty}$. The group $G_\omega$ is generated by the set of elements
$\mathcal{A}_{\omega}=\{a,b_{\omega},c_{\omega},d_{\omega} \}$ of order 2, with
$b_{\omega},c_{\omega},d_{\omega}$ commuting and generating the Klein 4-group (i.e.
$b_{\omega}c_{\omega}d_{\omega}=1$) (so  indeed the groups $G_\omega$ are 3-generated). For the definition of
these groups we address the reader to \cite{grigorch:degrees,grigorch:solved}. Originally $G_\omega$ were
defined similarly to $\mathcal{G}$ as groups acting on $[0,1]$ (with removed diadic rational points), but
more convenient language to work with them is via action on binary sequences (via identification of a point
from $[0,1]$ with its binary expansion), or via actions by automorphisms on a binary rooted tree $T$, when we
identify vertices of the tree with corresponding binary sequences.

 Let $Q$ be a subgroup of $\Xi$ generated by the elements $x=a,y=d,z=cac$. It is easy to check that $Q$ has index 2
 in $\Xi$ and has a presentation
 \[\langle x,y,z\|\ x^2,y^2,z^2,(xy)^4,(xz)^4,(yz)^4 \rangle.\]
 Therefore $Q$ is isomorphic to $\Phi$.

Let $\Omega_0 \subset \Omega$ be the subset consisting of sequences $\omega$ which contain each symbol
$0,1,2$ infinitely many times, $\Omega_1 \subset \Omega$ be the set of sequences which contain at least two
symbols from $\{0,1,2\}$ infinitely many times, and $\Omega_2=\Omega \setminus \Omega_1$ be the set of
sequences $\omega=\omega_1 \omega_2 \dots \omega_n \dots$ such that
$\omega_n=\omega_{n+1}=\omega_{n+2}=\dots$ starting with some coordinate $n$. Observe that all sets
$\Omega_0, \Omega_1, \Omega_2$ are invariant with respect to the shift $\tau$
\[\tau(\omega_1\omega_2\omega_3\dots) = \omega_2\omega_3\dots\]
 in the space of sequences. The
groups $G_{\omega}$ are virtually abelian for $\omega \in \Omega_2$, while the groups $G_{\omega}$, for
$\omega \in \Omega_1$ are just-infinite, branch groups of intermediate growth. Additionally, the groups
$G_{\omega}$, for $\omega \in \Omega_0$ are 2-groups. Profs of these facts are provided by
Theorems~2.1,~2.2,~8.1, and Corollary~3.2  in~\cite{grigorch:degrees}. One of important facts that will be used in the proof of the Theorem~\ref{coxeter} is that the set of growth degrees of groups $G_{\omega},
\omega \in \Omega_0$ has uncountable cardinality. The word problem for the family $G_{\omega}, \omega \in
\Omega_1$ can be solved by algorithm with oracle $\omega$ (i.e. the algorithm which uses the symbols of the sequence $\omega$ in its work), which we call branch algorithm because of its branching nature
\cite{grigorch:solved,grigorch:degrees}. Using this
 algorithm, or directly from the definition of groups $G_{\omega}$, it is easy to check that if $\omega$ begins with symbol $0$ then
 $(ad_{\omega})^4=1$, if $\omega=1\dots$, then $(ac_{\omega})^4=1$ and if $\omega=2\dots$, then
 $(ab_{\omega})^4=1$. As we can exclude any of $b_{\omega},c_{\omega},d_{\omega}$ from the generating set we see
 that each of the groups $G_{\omega}, \omega \in \Omega_1$ is a homomorphic image of $\Xi$. To simplify
 the situation we  assume that $\omega$ begins with $0$, so $(ad_{\omega})^4=1$. Let $\Omega_3\subset \Omega_0$ be the set
 of sequences which begin with symbol $0$. The proofs of results about growth in \cite{grigorch:degrees} allow to
 conclude that the set of growth degrees of groups from $\{G_{\omega}, \omega \in \Omega_3\}$ has uncountable
 cardinality. Moreover the same holds for any set of the form $w\Omega_0$, where $w $ is arbitrary finite binary
 sequence.

The results from \cite{grigorch:degrees} also show that the group $G_{\omega}, \omega \in \Omega_1$ is
abstractly
 commensurable with
$G_{\tau(\omega)} \times G_{\tau(\omega)}$ and therefore the growth of $G_{\omega}$ is equal to the square of
the growth of $G_{\tau(\omega)}$.

\section{Proof of the theorem}

\begin{proof} First we show that a Coxeter group $\mathcal{C}$, satisfying the condition of the
theorem~\ref{coxeter} can be mapped onto one of Coxeter groups $\Xi$, $\Phi$, $\Upsilon$, or $\Pi$. This will
reduce the proof to  these groups. Indeed everything will be deduced from the fact that the group $\Xi$
satisfies the conclusion of the theorem.
%The most crucial is the case of $\Xi$, because other cases
%can be reduced to it as we will see soon.

Assume that the graph $\mathcal{Z}$ is not a tree, so it contains a cycle of length $\geq 3$ consisting of
vertices $x_{i_1},x_{i_2},\dots, x_{i_k}$ for some $ 3 \leq k \leq n$. Taking the quotient of $\mathcal{C}$
by the normal subgroup generated by the generators $x_j$  which do not belong to this cycle, we can pass to
the case when the graph $\mathcal{Z}$ is a cycle. Taking the quotient by the relation $x_{i_1}=x_{i_2}$ (if
the length of the cycle is greater than 3) we make the cycle shorter. After finitely many steps of this type
we come to the case when the length of the cycle is 3. Then making the further factorization by replacing the
numbers $m_{i,j} \geq 8$ by $m_{i,j}=4$, we map $\mathcal{C}$ onto $\Phi$.

If $\mathcal{Z}$ is a tree, passing to an appropriate quotient reduces the situation to the case when the
graph $\mathcal{Z}$ looks like a ``segment''  (all vertices are of degree $\leq 2$) with 3 or 4 vertices, and
labeling of edges given by the set $\{4,8\}$ or $\{4,4,4\}$) respectively, or like a tripod ``Y'' (i.e. is a
tree with four vertices, one of degree 3 and three leaves) with all edges labeled by 4, which correspond to
the cases of groups $\Xi,\Upsilon$ and $\Pi$  respectively.

We already know from the previous section that  $\Xi$  has uncountably many quotients $G_{\omega}, \omega \in
\Omega_3$, with different types of growth which are branch just-infinite  2-groups. $\Phi$ is a subgroup of
index 2 in $\Xi$. Let $Q_{\omega}$ be the corresponding quotient of $\Phi$ in $G_{\omega}$  of index 2.
Obviously  $Q_{\omega}$  is 2-group and has the same growth type as $G_{\omega}$. The group $Q_{\omega}$ (as
well as $G_{\omega}$) acts on binary rooted tree $T$, and for branch property we need only to show that the
action is level transitive because for each $n$ the rigid stabilizer $rist_{Q_{\omega}}(n)$ has index $\leq
2$ in $rist_{G_{\omega}}(n)$. But $Q_{\omega}$ acts transitively on the first level and both projections of
$st_{Q_{\omega}}(1)$ are equal to the subgroup $R_{\omega}=\langle d_{\tau(\omega)}, ac_{\tau(\omega)}
\rangle $ which is of index 2 in $G_{\tau(\omega)}$. This subgroup also acts transitively on the first level
and both projections of $st_{R_{\omega}}(1)$ are equal to the group $G_{\tau^2(\omega)}$, which is branch.
Therefore $R_{\omega}$ and $Q_{\omega}$  act transitively by \cite[Theorem~4]{grigorch:jibranch} and are
branch groups. We conclude that $\Phi$ has uncountably many quotients satisfying conclusion of
theorem~\ref{coxeter}.

Now we are going to consider the case of $\Upsilon$. Let $\Lambda$ be a subgroup of index 2 in $\Phi$
generated by the elements $u=xy, v=xz$. Then $\Lambda$ has a presentation
\[\Lambda= \langle u,v \|\ u^4,v^4,(uv)^4 \rangle .\]
Consider the subgroup  $\mathcal{S}$ of  $\mathcal{G}$ generated by the elements $ad$ and $(ac)^2$. It is a
quotient of $\Lambda$ with respect to the map
\[  u \rightarrow ad, \qquad v \rightarrow (ac)^2 .\]

 Computations show that $\psi(st_{\mathcal{S}}(1))$ is a subgroup in $\mathcal{G} \times \mathcal{G}$
generated by the pairs $(b,b), (da,ad), (bac,da),(badac,(da)^2)$, and the projections of this subgroup on
each factor is the group $\langle b, ac\rangle =\langle b,ad\rangle$, which has index 2 in $\mathcal{G}$ and
is branch. Therefore $\mathcal{S}$ acts transitively on levels, is  branch, and just-infinite.

Let $S_{\omega}$, for $\omega \in \Omega_3$ be the subgroups of $G_{\omega}$ generated by $ad_{\omega}$ and
$(ac_{\omega})^2$. Then the relators of $\Lambda$ are also relators of  $S_{\omega}$ with respect to the map
\[  \nu :~x \rightarrow ad_{\omega},~v \rightarrow (ac_{\omega})^2. \]
The image  $\psi(st_{S_{\omega}}(1))$ is a branch subgroup $\langle b_{\omega}, ac_{\omega}\rangle$ of index
2 in $G_{\omega}$. Therefore $S_{\omega}$ is also branch, just-infinite and has the same growth type as
$G_{\omega}$.

 Let $\bar{\Lambda}$ be any of the 2-quotients $S_{\omega}$ of $\Lambda$ given by the
previous arguments, and let $u,v$ be the set of generators of $\bar{\Lambda}$ which are the images of the
generators of $\Lambda$ (we keep the same notation for them). Consider the group $\bar{\Lambda}_1$, acting on
binary rooted tree $T$, generated by the element $a$ of order two (permutation of two subtrees $T_0,T_1$ with
roots at the first level) and the elements $b=(1,v),c=(u,u), d=(v,v)$, where $u,v$ and the identity element
act on the left or right subtree respectively (in a same way they act on the whole tree; here we use the
self-similarity of the binary tree). Then $a$ commutes with $c$ and $d$, $b$ commutes with $d$,  and
$(ab)^4=(bc)^4=(cd)^4=1$, so the group is a quotient of $\Upsilon$. The $\psi$-image of stabilizer of the
first level of $\bar{\Lambda}_1$ is a subdirect product of $\bar{\Lambda} \times \bar{\Lambda}$ and contains
the group $D \times D$ where $D$ is the normal closure of $v$ in $\bar{\Lambda}$. As $\bar{\Lambda}$ is
just-infinite, $D$ has finite index in $\Lambda$. Therefore the growth of $\bar{\Lambda}_1$ is equal to the
square of the growth of $\bar{\Lambda}$.  It is clear that $\bar{\Lambda}_1$ is branch and just-infinite. As
the set of  squares of growth degrees $\bar{\Lambda}$ has uncountable cardinality, we are done with this
case.

Now consider the last case of the group $\Pi$. Let $G=G_{\omega}, \omega \in \Omega_3$ be a 2-group, whose
generators will be denoted, for simplicity, by $a,b,c,d$ instead of  $a,b_{\omega }, c_{\omega },d_{\omega
}$. Recall that $a$ acts by permutation of the two subtrees $T_0,T_1$ of the binary tree with roots on the
first level. Consider the group $V=\langle a, \bar{a}, \bar{b}, \bar{c} \rangle$, where $\bar{a}, \bar{b},
\bar{c}$ are automorphisms of the tree fixing the vertices of the first level whose $\psi$-images are $(a,1),
(1,b),(1,c)$ respectively (here again we  use the self-similarity of binary rooted tree identifying $T$ with
$T_0,T_1$). Then the generators $ a, \bar{a}, \bar{b}, \bar{c}$ are of order 2, $\bar{a}, \bar{b}, \bar{c}$
commute, and $(a\bar{a})^4=(a\bar{b})^4=(a\bar{c})^4=1$, so the group is a homomorphic image of the $\Pi$
with respect to the map
\[a\mapsto a,\quad b\mapsto\bar{b},\quad c\mapsto\bar{c},\quad
d\mapsto\bar{a}.\]

The $\psi$-image of $st_V(1)$ is a subdirect product of $G \times G$ and contains $A\times A$, where $A$ is
the normal closure of $a$ in $G$  ($A$ has finite index in $G$, as $G$ is just-infinite). $V$ acts
transitively on levels and therefore is branch and just-infinite. The growth of $V$ is the square of the
growth of $G_{\omega}$.  Therefore $\Pi$ has uncountably many quotients satisfying the statement of the
theorem.
\end{proof}

\section{Concluding remarks}

In 2006,  J.~Wiegold  raised the following  question in Kourovka Notebook~\cite[16.101]{kourovka06}. Do there
exist uncountably many infinite 2-groups that are quotients of the group
\[\Delta=\langle x,y \|\ x^2,y^4,(xy)^8 \rangle ?\]
The problem is motivated by the following comment by J.~Wiegold ``There certainly exists one, namely the
subgroup of finite index in Grigorchuk's first group generated by $b$ and $ad$; see (R.~I.~Grigorchuk, {\it
Functional Anal. Appl.}, {\bf 14} (1980), 41--43).''

Immediately after the appearance we informed one of the Editors of Kourovka Notebook, I.~Khukhro, that the
answer to the question is positive, and that the results of~\cite{grigorch:degrees} can be easily used to
provide a justification. Unfortunately, it took  some time for the author to write the corresponding text,
and he is finally presenting his arguments in this note. Different argument has been used recently in the
article \cite{minasyan_os:periodic} and the authors were notified of the approach given here (they
acknowledgment this fact at the end of Section 2).

Let $L$ be a subgroup of $\Xi$ generated by $x_1=ac,x_2=ad$. Then  $L$  is a subgroup of index 2 in $\Xi$,
has a presentation
\[L= \langle x_1,x_2 \|\ x_1^4,x_2^8,(x_2x_1^{-1})^2\rangle,\]
and therefore is isomorphic to the group $\Delta$ via the map $x \rightarrow x_1^{-1}x_2, y \rightarrow x_1$.
 Let $L_{\omega}$ be the subgroup of $G_\omega$ of index 2 generated by $ac_{\omega}$ and
$ad_{\omega}$. Then, if $\omega$ begins with $0$  (and so $(ad_{\omega})^4=1$), the group $L_{\omega}$ is a
homomorphic image of $\Delta$. As the set of growth degrees of  groups $L_{\omega}, \omega \in \Omega_0,
\omega=0w_2\dots$ has cardinality $2^{\aleph_0}$ we get the affirmative answer to the Wiegold question.
Obviously $L_{\omega}$ are 2-groups. One can show that they are branch and just-infinite as it is shown in
\cite{grigorch:wiegold}. Observe that alternatively the groups $L_{\omega}$  can be defined as groups
generated by elements $x=b_{\omega}, y=ad_{\omega}$ as was suggested  by Wiegold in the case of
$\mathcal{G}$, and that $L_{\omega}$ satisfy the relations
\[1=x^2=y^4=(xy)^8=(xy^2)^{16}\]
(the provided list of defining relations in not complete). It is unclear if the power 16 in the last relation
can be replaced by 8, i.e. if there is an infinite 2-generated 2-group the set of defining relations of which
starts with

\[1=x^2=y^4=(xy)^8=(xy^2)^{8}.\]

%in an 2-generated infinite 2-group, while keeping the powers 2,4,8 without change.

There are other approaches for construction of   infinite torsion quotients of Coxeter groups. For instance,
for those Coxeter groups which can be mapped onto non-elementary hyperbolic groups  (in Gromov sense
\cite{gromov:hyperbolic}), or which are ``large'' groups in the sense of S.~Pride
\cite{pride:height,edjevet_pr:largenes} (a group is ``large'' if it has a subgroup of finite index that can
be mapped onto a free group of rank 2),  the results and constructions from
\cite{ivanov_ol:97,olshanskii-o:large06,delzant_Grom:burnside,minasyan_os:periodic} can be used.

The criterion for  a Coxeter group defined  by a connected Coxeter graph to be non-elementary hyperbolic,
given by G.~Moussong in~\cite{moussong:phd}, requires that  each Coxeter subgroup generated by a subset
$\{x_i,x_j,x_k\}$ of three generators is a hyperbolic triangular group, i.e. a group isomorphic to the group
$T^{*}_{m,n,q}= \langle x,y,z \|\ x^2,y^2,z^2,(xy)^m=(xz)^n=(yz)^q\rangle$ with
\[\frac{1}{m}+ \frac{1}{n} + \frac{1}{q} < 1.\]

The groups $\Xi$ and $\Phi$ are non-elementary hyperbolic and, as was indicated by T.~Januszkiewicz, the
groups $\Upsilon$ and $\Pi$  can be mapped onto non-elementary hyperbolic groups. Therefore all these groups
have uncountably many homomorphic torsion images of bounded degree according to \cite{minasyan_os:periodic}.

Indeed, all Coxeter groups which are not virtually abelian are ``large", which is a particular case of the
results by G.~Margulis and E.~Vinberg from \cite{margulis_vin:large}. This fact was also proved independently
by C.~Gonciulea, as is indicated in the A.~Lubotzki's review [MR1748082 (2001h:22016)] to
\cite{margulis_vin:large}, but published only in a weaker form \cite{gonciulea:coxeter}. Therefore  in view
of the results from \cite{olshanskii-o:large06,minasyan_os:periodic}, for any prime number $p$ and any
Coxeter group $\mathcal{C}$ that is not virtually abelian, there is $2^{\aleph_0}$ pairwise non isomorphic
quotients of $\mathcal{C}$ which are residually finite virtually $p$-groups. It is pointed out by
T.~Januszkiewicz that it is possible that every Coxeter group that is not virtually abelian has a
non-elementary hyperbolic quotient (perhaps this is a known fact). If this is the case, then every Coxeter
group that is not virtually abelian has uncountably many torsion quotients of bounded exponent.

Finally let us formulate an open question. The $p$-groups ($p \geq 3$ is a prime) of Gupta-Sidki
\cite{gupta_s:burnside} are 2-generated, residually finite, branch, and just-infinite. Their generators $x,y$
satisfy the relations $x^p=y^p=(x^iy^j)^{p^2}=1, 1 \leq i,j \leq p-1$.
\begin{prob} Let $p \geq 5$ be a prime. Does there exists a residually finite p-group generated by two elements $x,y$ subject to
the relations $x^p=y^p=(x^iy^j)^{p}=1, 1 \leq i,j \leq p-1$? Can such a group   have additionally some other
finiteness properties (for instance in the spirit of theorem~\ref{coxeter})?
\end{prob}
The prime $p=3$ is excluded for obvious reasons.

Observe that the quotient $\mathcal{G}$ of $\Xi$ is a self-similar group (historically it is the first
example of a non elementary self-similar group; more on self-similar groups see in
\cite{bartholdi_gn:fractal,nekrash:self-similar}). The groups $Q=Q_{\varsigma}$, $ \bar{\Lambda}_1, V$ used
in the proof of theorem~\ref{coxeter}, which are quotients of groups $\Phi,\Upsilon$ and $\Pi$ respectively,
are not self-similar. It would be interesting to find self-similar torsion quotients of $\Phi,\Upsilon$ and
$\Pi$ if they exist (or to show that there is no such quotients).

%%%%%%%%%%%%%%%%%%%%%%%%%%%%%%%%%%%%%%%%%%%%

% ----------------------------------------------------------------
\def\cprime{$'$} \def\cprime{$'$} \def\cprime{$'$} \def\cprime{$'$}
  \def\cprime{$'$} \def\cprime{$'$} \def\cprime{$'$} \def\cprime{$'$}
  \def\cprime{$'$}
\providecommand{\bysame}{\leavevmode\hbox to3em{\hrulefill}\thinspace}
\providecommand{\MR}{\relax\ifhmode\unskip\space\fi MR }
\providecommand{\MRhref}[2]{%
  \href{http://www.ams.org/mathscinet-getitem?mr=#1}{#2}
} \providecommand{\href}[2]{#2}

\end{document}